\documentclass[12pt]{amsart}
\usepackage{amsmath,amssymb}
\newtheorem{theorem}{Theorem}
\newtheorem{proposition}[theorem]{Proposition}
\newtheorem{lemma}[theorem]{Lemma}
\newtheorem{corollary}[theorem]{Corollary}
\newtheorem{definition}[theorem]{Definition}
\newenvironment{proof*}{\vskip 2mm\noindent {}}{\hfill $\Box$ \vskip 2mm}
\def\R{\mathbb R}
\def\C{\mathbb C}
\def\D{\mathbb D}

\def\eps{\varepsilon}
\renewcommand{\Im}{\operatorname{\mbox{Im}}}
\renewcommand{\Re}{\operatorname{\mbox{Re}}}
\newcommand{\diam}{\operatorname{\mbox{diam}}}

\title[Comparison of the real and the complex Green functions]{Comparison of the real and
the complex Green functions, and sharp estimates of the Kobayashi distance}

\author{Nikolai Nikolov and Pascal J. Thomas}

\address{N.~Nikolov\\ Institute of Mathematics and Informatics\\Bulgarian Academy
of Sciences\\Acad. G. Bonchev 8, 1113 Sofia, Bulgaria\newline
\indent Faculty of Information Sciences\\
State University of Library Studies and Information Technologies\\
Shipchenski prohod 69A, 1574 Sofia, Bulgaria}\email{nik@math.bas.bg}

\address{P.J.~Thomas\\
Institut de Mathématiques de Toulouse ; UMR5219 \\
Université de Toulouse ; CNRS \\
UPS IMT, 31062 Toulouse Cedex 9, France} \email{pascal.thomas@math.univ-toulouse.fr}



\thanks{This paper was started while the first-named author was an invited professor at
the Paul Saba\-tier University, Toulouse in May-June 2016.}

\begin{document}

\begin{abstract}{}
We extend the upper estimates obtained by M. Carlehed \cite{Ca} and B.-Y. Chen \cite{Ch}
about the ratio of the classical and pluricomplex Green functions to the case of
$\mathcal C^2$-smooth locally $\C$-convexifiable
domains of finite type. We also give some lower estimates. In order to obtain those
results, and because it is of independent interest, we refine and unify some classical
estimates about the Kobayashi distance and the Lempert function in such domains.
\end{abstract}

\maketitle

\section{Introduction and results}

\subsection{Green functions.}
Two kinds of Green functions can be defined on a domain $D \subset \C^n \cong \R^{2n}$, $n\ge 2$: the usual one,
related to harmonic (or subharmonic) functions when $D$ is seen as subdomain of $\R^{2n}$, and the pluricomplex
Green function (see e.g. \cite{JP}), related to plurisubharmonic functions.

The pluricomplex Green function
originated with the work of Lempert \cite{Le}, Klimek \cite{Kli}, Lelong \cite{Lel},
among others, and is the subject
of many recent works, see for instance \cite{He}, \cite{RT}.

Let $G_D$ stand for the usual \emph{Green function} at a pole $w$ in $D \subset \R^m$, $m\ge 3$, given by
$$
G_D (z,w) = \sup\left\{ u(z): u \in SH_-(D),\ u = |\cdot-w|^{-m+2} + O(1)\right\}.
$$
Let $g_D$ stand for the \emph{pluricomplex Green function} at a pole $w$ in $D \subset \C^n$, $n\ge 2$, given by
$$
g_D (z,w) = \sup\left\{ u(z): u\in PSH_-(D),\ u = \log|\cdot-w| + O(1)\right\}.
$$
Here $SH_-(D)$ and $PSH_-(D)$ stand for negative subharmonic, resp. plurisubharmonic, functions on $D$.

Note that for $n=1$ the second extremal problem also gives the usual Green functions for the Laplacian on $\R^2$.

The respective behavior of those two functions were compared by M. Carlehed \cite{Ca} and B.-Y. Chen \cite{Ch}.
In the present paper, we extend their results to a wider class of domains, and
give some improved estimates
for various holomorphic invariants such as the Kobayashi distance in that class of domains.

We would like to thank the referee for his very careful reading of our manuscript and
several useful suggestions. 

\subsection{Domains in $\C^n$.}
In order to state the results, we need to define some geometric properties of a domain in $\C^n$.
From now on, we assume that $n\ge 2.$ As usual, we say that $\partial D$, or $D$, is $\mathcal C^k$-smooth if $D=\{ \rho <0\}$,
where $\rho$ is a defining function of class $\mathcal C^k$ on $\overline D$
such that $\nabla \rho$ does not vanish on $\partial D$. A $\mathcal C^2$-smooth domain
is \emph{strictly pseudoconvex} if the complex Hessian of $\rho$ restricted to the complex tangent
space at every point of $\partial D$ is positive definite.

A domain $D$ is \emph{$\C$-convex} if any non-empty intersection of $D$ with a complex line
is connected and simply connected. If $D$ is bounded and $\mathcal C^1$-smooth, this is equivalent
to being \emph{lineally convex}, that is to say,
for any $z\notin D$, there exists a complex hyperplane $H$ through $z$ such that
$D\cap H =\emptyset$. For more on those two notions, see e.g. \cite{APS}.

A domain $D$ is \emph{$\C$-convexifiable} if $D$ is biholomorphic to a $\C$-convex domain.

A domain $D$ is \emph{locally ($\C$-)convexifiable}, if for any $a\in\partial D$,
there exist a neighborhood $U$ of $a$ and a holomorphic embedding $\Phi:U\to\C^n$
such that $\Phi(D\cap U)$ is a ($\C$-)convex, domain.

It is well-known that any strictly pseudoconvex domain is locally convexifiable.

The \emph{type} of a smooth boundary point $a$ of a domain $D$ is the
supremum over the orders of contact of the one-dimensional analytic varieties through $a$ with $\partial D$
(possibly $\infty$). The type of a smooth domain $D$ is defined as the supremum over the types of all boundary
points of $D.$ For instance, the bounded domains of type $2$ are exactly the strictly pseudoconvex domains.
Also, the types of the pseudoconvex domains are even numbers or $\infty.$
If the domain is $\C$-convex, the type does not change, considering complex lines instead of
varieties (see e.g. \cite[Proposition 6]{NPZ}).

\subsection{Notations and auxiliary quantities.}

We will systematically use the following notations : $A \gtrsim B$ means that there
is a constant $C>0$ such that $A \ge C B$;
$A\asymp B$ means that $A \gtrsim B$ and $B \gtrsim A$;
and $A \sim B$ means that $A/B \to 1$. What the constants depend on, and
in which sense the limit is taken, will be made clear from context.

The Green functions we consider take negative values and, when $\partial D$ is smooth enough, tend to $0$
at the boundary.  A typical negative plurisubharmonic function is $\log |f|$, where $f$ is a
holomorphic function bounded by $1$; so it will be convenient to consider $e^{g_D}$.  Consideration
of the Poincar\'e distance $p$ in the unit disc $\D$,
$p(w,z)= \tanh^{-1} \left|\frac{z-w}{1-\bar z w} \right|$, makes it expedient to
consider $\tanh^{-1} e^{g_D}$.

We give a unified convention.

\begin{definition}
\label{convention}
Given any continuous function $f : D \to (-\infty, 0)$, we write
\begin{equation}
\label{star}
f^*:= e^f, \mbox{ so } f^* : D \to (0, 1),
\end{equation}
\begin{equation}
\label{tilde}
\tilde f := \tanh^{-1} f^* =\tanh^{-1} e^{f} = \frac12 \log \frac{1+e^{f}}{1-e^{f}},
\mbox{ so } \tilde f : D \to (0, \infty).
\end{equation}
\end{definition}
Conversely, $f^*= \tanh \tilde f = \frac{e^{2\tilde f}-1}{e^{2\tilde f}+1}$,
and $f= \log f^*$.

Elementary calculations give:
\begin{lemma}
\label{tech}
\
\begin{enumerate}
\item[(i)]
Suppose that $f \to 0^-$, or equivalently $f^*\to 1^-$, or equivalently $\tilde f \to \infty$. Then
$1-f^* \sim -f,$ $\tilde f \sim -\frac12\log(-f),$ and
$f \sim -2 e^{-2 \tilde f }$; in particular if $\tilde f= \log t$, then $f \sim -\frac2{t^2 }$.

\item[(ii)]
Suppose that $f \to -\infty$, or equivalently $f^*\to 0^+$, or equivalently $\tilde f \to 0^+$.
Then $\tilde f \sim f^*$ and $f = \log \tilde f + O(1)$.
\end{enumerate}
\end{lemma}

\subsection{The ratio of the Green functions.}
Our first main result is the extension to the case of locally $\C$-convexifiable domains of
a theorem proved in the case of locally convexifiable domains \cite[Theorem 1]{Ch}.

\begin{theorem}
\label{greenratio}
Let $D \subset \C^n$ be a bounded, smooth, locally $\C$-convexi\-fiable
domain of type $2m.$ Then there exists $C>0$ such that
$$
\frac{g_D(z,w)}{G_D(z,w)} \le C |z-w|^{2(n-2m)},\quad z,w\in D,\ z\neq w.
$$
\end{theorem}

For $z \in D$, let $\delta_D(z) := \min \left\{ |z-w|: w \notin D \right\}$
(the distance to the boundary). Any bounded, $\mathcal C^{1,1}$-smooth domain $D$
is of \emph{positive reach}, that is to say, there exists $\delta_0 >0$ such that for any
$z\in D$ with $\delta_D(z)< \delta_0$, there exists a unique
point $\pi(z) \in \partial D$ such that $|z-\pi(z)|=\delta_D(z)$.

Recall the following estimate of $G_D,$ when $D$ is bounded, $\mathcal C^{1,1}$-smooth domain
in $\R^m,$ $m\ge 3$ (see e.g. \cite[(7)]{Sw}):
\begin{equation}\label{sw}
c_1 G_D(z,w)\le-\min\Bigl\{\frac{1}{|z-w|^{m-2}},\frac{\delta_D(z)\delta_D(w)}{|z-w|^m}\Bigr\}
\le c_2 G_D(z,w),
\end{equation}
where $c_1,c_2>0$ are constants, and $z,w\in D.$

The proof of Theorem \ref{greenratio} will rely on the second inequality
in \eqref{sw}, and the following precise estimate of the pluricomplex Green function $g_D$ which
is sensitive in both extreme cases: $g_D\to 0$ and $g_D\to-\infty.$

\begin{theorem}
\label{greengrowth}
Let $D$ be as in Theorem \ref{greenratio}. Then there exists $C>0$ such that for any $z,w\in D$,
\begin{equation}\label{growth}
\tilde g_D (z,w)\ge
m\log \left( 1+C \frac{|z-w|}{\delta_D(z)^{1/2m}} \right)\left( 1+C \frac{|z-w|}{\delta_D(w)^{1/2m}} \right).
\end{equation}
\end{theorem}

In the more general case of a bounded, smooth, pseudoconvex domain of finite type, a weaker
estimate is proved by G. Herbort \cite[Theorem 1.1]{He}.

The proof of Theorem \ref{greengrowth} will be based on the respective local estimates, covering the cases where
either the pole or the argument tends to a boundary point.

\begin{theorem}\label{greenlocal}
Let $D\subset \C^n$ be a bounded domain, which is smooth and locally $\C$-convexifiable
near point $a \in \partial D$ of type $2m$.
Then there exist a neighborhood $U$ of $a$ and $C>0$ such that
\begin{equation}\label{pole}
\tilde g_D (z,w) \ge m \log \left( 1+C \frac{|z-w|}{\delta_D(w)^{1/2m}} \right),\quad z \in D,\ w \in D\cap U,
\end{equation}
\begin{equation}\label{argu}
\tilde g_D (z,w) \ge m \log \left( 1+C \frac{|z-w|}{\delta_D(z)^{1/2m}} \right), \quad z \in D\cap U,\ w \in D.
\end{equation}
\end{theorem}
In the particular case when $D$ is locally convexifiable, 
similar but weaker estimates than those in the above two theorems are contained in \cite{Ch}.

\subsection{Other holomorphic invariants.}
We will use other holomorphically contractive functions, with notations
sometimes slightly different from those of the standard reference \cite{JP},
to stay in line with the convention from Definition \ref{convention}. In particular,
note that the Kobayashi pseudo-distance in a domain $D$ will be called $\tilde k_D$,
while $k_D:= \log \tanh \tilde k_D \in (-\infty,0)$.  This is because our main focus
is on (negative-valued) Green functions.

Let $D \subset \C^n$, and $z,w \in D$.

The \emph{Lempert function} is given by
$$
\tilde l_D (z,w) := \inf \left\{ p(\zeta, \omega) : \zeta, \omega \in \mathbb D,
\exists \varphi \in \mathcal O (\D, D): \varphi (\zeta) = z, \varphi (\omega) = w
\right\}.
$$
With the notation convention from Definition \ref{convention}, this means that
$$
l_D^* (z,w) := \inf \left\{ \left|\frac{\zeta- \omega}{1-\bar \zeta \omega} \right| : \zeta, \omega \in \mathbb D,
\exists \varphi \in \mathcal O (\D, D): \varphi (\zeta) = z, \varphi (\omega) = w
\right\},
$$
and that $l_D (z,w) =\log l_D^* (z,w) \in (-\infty,0)$, a quantity that is easier to
compare with the Green function.

The \emph{Kobayashi-Royden (pseudo)metric} applied to a vector $X \in \C^n$ is given by
$$
\kappa_D (z;X) :=  \inf \left\{ \lambda >0  :
\exists \varphi \in \mathcal O (\D, D): \varphi (0) = z, \lambda \varphi' (0) = X
\right\}.
$$

The  \emph{Kobayashi (pseudo)distance} is the largest pseudodistance dominated by the
Lempert function. It is also given
by
$$
\tilde k_D (z,w) := \inf_\gamma\int_0^1 \kappa_D (\gamma(t);\gamma'(t)) dt,
$$
where the infimum is taken over all $\mathcal C^1$-smooth curves $\gamma:[0,1]\to D$
with $\gamma(0)=z$ and $\gamma(1)=w.$  Then $k_D(z,w)= \log \tanh \left( \tilde k_D (z,w)\right)$.

We have that
\begin{equation}\label{ineq} k_D \le l_D,\quad g_D\le l_D.
\end{equation}
Lempert's celebrated theorem \cite{Le} implies that in the case of a convex domain, those
are all equalities. This extends to the case of bounded, $\mathcal C^2$-smooth, $\C$-convex
domains \cite{Ja}. No inequality holds in general between $\tilde k_D$ and $\tilde g_D$;
and while $\tilde k_D$ is symmetric in its arguments, $\tilde g_D$ is not always so, but
we will see that under our hypotheses, they exhibit similar behavior.

\subsection{Lower estimates of the Kobayashi distance}

\begin{theorem}
\label{kobagrowth}
Let $D$ be as in Theorem \ref{greenratio}. Then there exists $C>0$ such that for
any $z,w \in D,$
\begin{equation}
\label{kobag}\tilde k_D (z,w) \ge
m\log\left( 1+C \frac{|z-w|}{\delta_D(z)^{1/2m}} \right)
\left( 1+C \frac{|z-w|}{\delta_D(w)^{1/2m}} \right).
\end{equation}
\end{theorem}

This will follow from the corresponding local sharp result.

\begin{theorem}
\label{kobalocal}
Let $D\subset \C^n$ be a domain, which is smooth and locally $\C$-convexifiable
near a point $a \in \partial D$ of type $2m$.  Then there exist a neighborhood $U$ of $a$ and
$C>0$ such that for any $z \in D\cap U$, $w \in D$,
\begin{equation}
\label{kobaloc}
\tilde k_D (z,w) \ge m \log \left( 1+C \frac{|z-w|}{\delta_D(z)^{1/2m}} \right).
\end{equation}
\end{theorem}

\subsection{Upper bounds for the Lempert function and sharpness of the results.}
The next propositions
(inspired by the examples in \cite[p.~404]{Ca} and \cite[p.~35]{Ch}) and
\eqref{ineq} show that the exponents in all the above theorems are optimal.

\begin{proposition}
\label{sharp}
Let $D\subset \C^n$ be a domain, which is smooth and $\C$-convex
near a point $a \in \partial D$ of type $2m$.  Denote by $n_a$ the inner normal half-line
to $\partial D$ at $a.$ If $a$ is of type $2m,$ there exist a unit vector $X\in T_a^\C\partial D$
and $C>0$ such that for all $z\in n_a$, close enough to $a$, and
all $w \in D$ such that $\frac{z-w}{|z-w|}=X$ and $C\frac{|z-w|}{\delta_D(z)^{1/2m}}<1$, then
$$
l_D^*(z,w)
\le
C\frac{|z-w|}{\delta_D(z)^{1/2m}}.
$$
If $a$ is of infinite type, the last inequality holds for any $m\in\mathbb N$ with $C=C_m.$
\end{proposition}

We then have the following result characterizing the type of a point.

\begin{corollary}
Let $D\subset \C^n$ be a domain, which is smooth and locally $\C$-convexifiable
near a point $a \in \partial D.$ Then any of the inequalities \eqref{pole}, \eqref{argu}
and \eqref{kobaloc} holds if and only if $a$ is of type at most $2m.$
\end{corollary}

The next results are related to the converse of Theorem \ref{greengrowth}.

\begin{proposition}
\label{seq} Let $D \subset \C^n$ be a bounded, smooth, locally $\C$-convexi\-fiable
domain. If $D$ is of type $2m,$ there exist sequences $(z^j),(w^j)\subset D$ and $c>0$ such that
$|z^j-w^j|\rightarrow 0$ and
$$
\frac{g_D(z^j,w^j)}{G_D(z^j,w^j)}\ge c|z^j-w^j|^{2(n-2m)},\quad j\in\mathbb N.
$$
If $D$ is of infinite type, the last inequality holds for any $m\in\mathbb N$ with
$(z^j)$, $(w^j)$ and $c$ depending on $m.$
\end{proposition}

Theorem \ref{greenratio} and Proposition \ref{seq} imply the following characterizations
of the type of a domain.

\begin{corollary}
\label{type}
Let $D \subset \C^n$ be a bounded, smooth, locally $\C$-convexi\-fiable
domain. Then:
\smallskip

\noindent(i) there exists $C>0$ such that
$$
\frac{g_D(z,w)}{G_D(z,w)} \le C |z-w|^{2(n-2m)},\quad z,w\in D,\ z\neq w.
$$
if only if $D$ is of type at most $2m;$
\smallskip

\noindent(ii) the ratio $g_D/G_D$ is bounded from above if and only if
$D$ is of type at most $n.$
\end{corollary}

If $m=1,$ the condition about $\C$-convexity is superfluous.

\begin{proposition}\label{strict}
Let $D \subset \C^n$ be a bounded, $\mathcal C^2$-smooth
domain. Then there exists $C>0$ such that
\begin{equation}
\label{str}
\frac{g_D(z,w)}{G_D(z,w)} \le C |z-w|^{2n-4},\quad z,w\in D,\ z\neq w.
\end{equation}
if and only if $D$ is strictly pseudoconvex.
\end{proposition}

In dimension $2,$ this proposition says that the ratio $g_D/G_D$ is bounded from above
if only if $D$ is strictly pseudoconvex. By Corollary \ref{type}, this is not true if $n\ge 4.$

\begin{proposition}\label{three}
Let $D \subset \C^3$ be a bounded, $\mathcal C^3$-smooth domain. Then the ratio $g_D/G_D$ is
bounded from above if only if $D$ is strictly pseudoconvex.
\end{proposition}

\smallskip

It is natural to ask which upper bounds can be given for the functions $g_D$ and $k_D$,
and indeed, many results for $k_D$ have been given in that direction, see for instance
\cite{NA}. To get estimates from above, using \eqref{ineq}, it will be enough to bound $\tilde l_D (z,w)$.

\begin{proposition}
\label{upper}
Let $D \subset \C^n$ be a bounded, $\mathcal C^2$-smooth, locally $\C$-convexi\-fiable
domain. Then there exists $C>0$ such that
\begin{equation}
\label{upperbd}
\tilde l_D (z,w)
\le \log \left( 1+C \frac{|z-w|}{\delta_D(z)^{1/2}\delta_D(w)^{1/2}} \right),\quad z,w\in D.
\end{equation}
\end{proposition}

This proposition shows that the factor $m$ in Theorems \ref{greengrowth}--\ref{kobalocal}
is sharp. On the other hand, these theorems show that the exponent $1/2$ in
Proposition \ref{upper} is optimal.

Proposition \ref{upper}, \eqref{sw}, \eqref{ineq}, and Lemma \ref{tech} also imply the following:
\begin{corollary}
Let $D$ be as in Proposition \ref{upper}.
Then there exists $C>0$ such that
\begin{equation}
\label{uppr}
\frac{g_D(z,w)}{G_D(z,w)} \ge C |z-w|^{2n-2},\quad z,w\in D,\ z\neq w.
\end{equation}
\end{corollary}

We already know from \cite[Theorem 2]{NPT1}
that if $D$ is a bounded, $\mathcal C^{1+\varepsilon}$ domain in $\C^n,$
then a weaker estimate than \eqref{upperbd} holds:
\begin{equation}\label{lem}
\tilde l_D (z,w)
\le \log \frac{C}{\delta_D(z)^{1/2}\delta_D(w)^{1/2}}.
\end{equation}
It would be interesting to know if \eqref{upperbd} and, hence, \eqref{uppr}
remain true in this general case.

\

The rest of the paper is organized as follows: Section \ref{proofupper}
contains the proofs of Propositions \ref{sharp}, \ref{seq}, \ref{strict},
\ref{three}, and \ref{upper}, Section \ref{estkoba} -- the proofs of Theorems
\ref{kobagrowth} and \ref{kobalocal}, Section \ref{estgreen} -- the proof
of Theorem \ref{greenlocal}, and Section \ref{proofgreen} --
the proofs of Theorem \ref{greenratio} and \ref{greengrowth}.

\section{Proofs of Propositions \ref{sharp}, \ref{seq}, \ref{strict}, \ref{three}, and \ref{upper}}
\label{proofupper}

\noindent{\it Proof of Proposition \ref{sharp}.}
By \cite[Propositions 4 and 6]{NPZ}, if $a$ is of type at least $2m,$
there exist a neighborhood $U$ of $a,$ a unit vector $X\in T_a^\C\partial D,$ and $C>0$
such that the distance $\delta_D(z;X)$ from $z\in D\cap U\cap n_a$ to $\partial D$ in direction $X$
verifies $\delta_D(z;X) \ge C\delta_D(z)^{1/2m}.$ If $a$ is of infinite type, the last holds for any
$m\in\mathbb N$ with $C=C_m.$ Let $D_{z,X}:= \{z+tX: |t|< \delta_D(z,X)\}$.
Let $w-z = \lambda X$.
It remains to observe that if $r>1$ and $r|\lambda|<\delta_D(z;X),$ then,
recalling that $0<l_D^*<1$ with the notations from Definition \ref{convention},
$$
 l_D^* (z,z+\lambda X)\le l_{D_{z,X}}^* (z,z+\lambda X) =
 \frac{|\lambda|}{\delta_D(z;X)} \le \frac1{r} <1. \quad \qed
$$

\begin{proof*}{\it Proof of Proposition \ref{seq}.} Let $D$ be of type $2m.$
Choose a point  $a\in\partial D$ of type $2m.$
There exist a neighborhood $U_0$ of $a$ and a holomorphic embedding $\Phi:U_0\to\C^n$
such that $\Omega:=\Phi(D\cap U_0)$ is a $\C$-convex domain. Set $u'=\Phi(u).$ Since
$|z'-w'|\asymp|z-w|$ and $\delta_\Omega(u')\asymp\delta_D(u)$ for $u,z,w\in U_1\Subset U$, and
$l_\Omega(z',w')\ge l_D(z,w),$ we may assume that $D$ is $\C$-convex.

Let $X$ be as in Proposition \ref{sharp}. Using e.g. a smooth defining function of
$D$ near $a,$ one may find a neighborhood $U$ of $a$ and $C>1$ such that
if $z\in D\cap U\cap n_a$ and $w=z+\lambda X,$ $C|\lambda|<\delta_D(z)^{1/2m},$ then
$\delta_D(z)=|z-a|<C\delta_D(w).$ Changing $U$ and $C$ (if necessary), we may apply Proposition
\ref{sharp} to find sequences $(z_j),(w_j)\to a$ such that $\delta_D(z_j)\asymp\delta_D(w_j)\asymp
|z_j-w_j|^{2m}$ and $\tilde l_D(z_j,w_j)\lesssim 1.$

This and the inequalities \eqref{sw} and \eqref{ineq} imply the desired result in the finite type case.

Let $D$ be of infinite type. Since $D$ is locally $\C$-convexifiable, there exists a point
$a\in\partial D$ of infinite type. Then, for any $m\in\mathbb N,$ we may proceed as above.
\end{proof*}

\noindent{\it Proof of Proposition \ref{strict}.} Strict pseudoconvexity implies local
convexifiability and, hence, \eqref{str} by Theorem \ref{greenratio}.

To prove the converse, we will proceed similarly to the proof of Proposition \ref{seq}.

Assume that the ratio $g_D/G_D$ is bounded from above, and
$a\in\partial D$ is not a strictly pseudoconvex point.

After an affine change of coordinates, we may suppose that $a=0$
and that $D$ is defined near $0$ by
$$
\Re(z_1+c_1z_2^2)\\
+c_2|z_2|^2+o(|\Im(z_1)|+|z_2|^2+|z''|)<0$$
where $c_2\le 0$.

It follows by \eqref{sw} that $g_D(z,w_0)\to 0$ as $z\to\partial D$
and hence $D$ is a pseudoconvex domain. This implies that $c_2=0.$

Let $\Phi(z)=(z_1+c_1z_2^2,z_2,z'').$ Then $G:=\Phi(D)$ is given near $0$ by
$$
\Re z_1+o(|\Im(z_1)|+|z_2|^2+|z''|)<0.
$$
Now it is easy to find sequences $\R_- \times \{0'\}\supset(z^j)\to 0$
and $(\lambda_j)\to\infty$ such that $G\ni w^j=z^j+\lambda_j\delta_G(z^j)^{1/2} e_2,$
and $2|z^j-w^j|<\delta_G(z^j;e_2)$, where $e_2:=(0,1,0, \dots, 0)$.

Because the order of contact of $\partial G$  and $\C e_2$ at $0$ is at least $2$, $|\delta_G (z^j)-\delta_G (w^j)|=O(|z^j-w^j|^2),$ so
$$
\frac{\delta_G(z^j)\delta_G(w^j)}{|z^j-w^j|^4}
\lesssim
\frac{\delta_G(z^j)^2}{|z^j-w^j|^4} +  \frac{\delta_G(z^j)}{|z^j-w^j|^2}
\to 0 \mbox{ and } l^*_G(z^j,w^j)<\frac12.
$$


If $\tilde z^j=\Phi^{-1}(z^j)$ and $\tilde w^j=\Phi^{-1}(w^j),$
then the inequalities $g_D\le l_D\le l_{D\cap U}$
and \eqref{sw} easily lead to the contradiction
$$\frac{g_D(\tilde z^j,\tilde w^j)}{G_D(\tilde z^j,\tilde w^j)}
|\tilde z^j-\tilde w^j|^{4-2n}\to\infty.\quad\qed
$$

\noindent{\it Proof of Proposition \ref{three}.}
As above, strict pseudoconvexity implies that
$$\frac{g_D(z,w)}{G_D(z,w)}\lesssim|z-w|^2\lesssim 1,\quad z,w\in D,\ z\neq w.$$

For the converse, assume that the ratio $g_D/G_D$ is bounded from above, and
$a\in\partial D$ is not a strictly pseudoconvex point.

After biholomorphic changes of variables similar to that in the proof of Proposition \ref{strict},
we may suppose that $D$ is defined near $a=0$ by
$$
\Re(z_1+c_3z_2^3+c_4z_2^2\overline{z_2})+o(|\Im(z_1)|+|z_2|^3+|z_3|)<0,
$$

Again by pseudoconvexity, $c_4=0.$ Let $\Psi(z)=(z_1+c_3z_2^3,z_2,z_3)$ and
Then $E:=\Psi(D)$ is defined near $0$ by
$$
\Re(z_1)+o(|\Im(z_1)|+|z_2|^3+|z_3|)<0.
$$

We may proceed as at the end of the proof of Proposition \ref{strict} to get a contradiction,
finding sequences $(z^j),(w^j)\to 0$ and $(\lambda^j)\to\infty$ such that
$w^j=z^j+\lambda^j\delta_E(z^j)^{1/3}e_2$, $l^*_E(z^j,w^j)<\frac12$,
and
since the order of contact of $\partial E$ at $0$ and $\C e_2$ is at least $3$,
$|\delta_E (z^j)-\delta_E (w^j)|=O(|z^j-w^j|^3),$ so
$$\frac{\delta_E(z^j)\delta_E(w^j)}{|z^j-w^j|^6}\to 0 \mbox{ and }
\frac{g_D(z^j,w^j)}{G_D(z^j,w^j)}\to\infty.\quad\qed$$




\noindent{\it Proof of Proposition \ref{upper}.}
By \eqref{lem}, for a given $\eps_0>0$, \eqref{upperbd} follows for $|z-w|\ge \eps_0$.
If $\min\left(\delta_D(w),\delta_D(z)\right) \ge \eps_0$, \eqref{upperbd} also follows, trivially.
So we may assume, by symmetry of the function, that $\delta_D(z) \le \delta_D(w) \le 2\eps_0$.

For any $a \in \partial D$, we may choose a bounded neighborhood $U_0$ of $a$
such that $D\cap U_0$ is $\C$-convexifiable and $\mathcal C^2$-smooth
(see \cite[Proposition 3.3]{NPT2}),
and that the projection $\pi$ to $\partial D$ is well defined on $U_0$.
Choose neighborhoods of $a$, $U_2 \Subset U_1$, such that $D\cap U_1 \Subset D\cap U_0$,
and $\eps_1 >0$ such that $z \in D\cap U_1$ and $\delta_D(z) \le \eps_1$ imply
$\delta_{D\cap U_0}(z) = \delta_D(z)$. We can cover $\partial D$ by a finite collection
of the $U_2$, and choose $\eps_0>0$ so that for any $z,w$ such that
$\delta_D(z) \le \delta_D(w) \le 2\eps_0$ and $|z-w|\le \eps_0$,
then $z \in U_2$, $w \in U_1$ (for some $a \in \partial D$) and
$\delta_{D\cap U_0}(z) = \delta_D(z)$, $\delta_{D\cap U_0}(w) = \delta_D(w)$.

Given $z,w$ as above, $\tilde l_D(z,w) \le \tilde l_{D\cap U_0}(z,w)$.

Then, by Lempert's Theorem, $\tilde l_{D\cap U_0}=\tilde k_{D\cap U_0}$, and by \cite[Corollary 8]{NA},
\begin{eqnarray*}
\tilde k_{D\cap U_0}(z,w)&\le&
\log \left( 1+C \frac{|z-w|}{\delta_{D\cap U_0}(z)^{1/2}\delta_{D\cap U_0}(w)^{1/2}} \right)\\
&=&\log \left( 1+C \frac{|z-w|}{\delta_D(z)^{1/2}\delta_D(w)^{1/2}} \right).\qed\quad
\end{eqnarray*}

\section{Proofs of Theorems \ref{kobagrowth} and \ref{kobalocal}}
\label{estkoba}
\begin{proof*}{\it Proof of Theorem \ref{kobagrowth}.}
Under the hypotheses of Theorem \ref{kobagrowth}, Theorem \ref{kobalocal}
and an compactness argument show that there is $\delta_0>0$ such that
\eqref{kobaloc} holds uniformly for $z,w\in D$ if $\delta_D(z)<2\delta_0.$
By symmetry, it is enough to consider three cases.

{\bf Case 1.} $\delta_D(z)\ge\delta_0,$ $\delta_D(w)\ge\delta_0.$

Then \eqref{kobag} follows from the inequality $\tilde k_D(z,w)\gtrsim |z-w|$,
valid on any bounded domain.

{\bf Case 2.} $\delta_D(z)<\delta_0,$ $\delta_D(w)\ge2\delta_0.$

Then $ \frac{|z-w|}{\delta_D(z)^{1/2m}}\gtrsim 1\gtrsim\frac{|z-w|}{\delta_D(w)^{1/2m}}$
and \eqref{kobag} follows by \eqref{kobaloc} (with bigger $C$).

{\bf Case 3.} $\delta_D(z)<\delta_0,$ $\delta_D(w)<2\delta_0.$

For any $\eps >0$, choose a curve
$\gamma$ so that its Kobayashi-Royden length is bounded by
$(1+\eps)\tilde k_D (z,w)$.
Choose a point $u\in\gamma$ such that
$|z-u|=|u-w| \ge \frac12 |z-w|.$
Then the definition of the Kobayashi distance and
\eqref{kobaloc} applied to $ (z,u)$ and $ (w,u)$ imply
\begin{multline*}
(1+\eps)\tilde k_D (z,w) \ge \tilde k_D (z,u)+\tilde k_D (u,w)
\\
\ge  m \log \left( 1+C \frac{|z-w|}{2 \delta_D(z)^{1/2m}} \right) +
m \log \left( 1+C \frac{|z-w|}{2 \delta_D(w)^{1/2m}} \right),
\end{multline*}
which, replacing $C$ by $C/2$, finishes the proof.
\end{proof*}

\begin{proof*}{\it Proof of Theorem \ref{kobalocal}.}
There exist a neighborhood $U_0$ of $a$ and a holomorphic embedding $\Phi:U_0\to\C^n$
such that $\Omega:=\Phi(D\cap U_0)$ is a $\C$-convex domain. Let $U_1$ and $U_2$ be neighborhoods
of $a$ such that $U_1\Subset U_2\Subset U_0.$ Let $z\in D\cap U_1.$

{\bf Case 1. $|z-w|^{2m}\le\delta_D(z)$.}

Since $\log(1+x)\le x,$ it is enough to prove that
\begin{equation}
\label{mest}
\tilde k_D (z,w) \gtrsim \frac{|z-w|}{\delta_D(z)^{1/2m}}.
\end{equation}
Let $\tilde k_D (D\cap U_1, D\setminus U_2)=:C_1>0.$
We may assume that $\tilde k_D (z,w)<C_1.$ Then
a curve connecting $z$ and $w$  of Kobayashi-Royden length $<C_1$
must lie inside $U_2.$
Since
$$\kappa_D(u,X)\gtrsim\kappa_{D\cap U_0}(u,X),\quad u\in U_2, X\in \C^n$$
(see e.g. \cite[Proposition 7.2.9]{JP}),
then $\tilde k_D(z,w)\gtrsim \tilde k_{D\cap U_0}(z,w)$.

From now on, we estimate $\tilde k_{D\cap U_0}(z,w).$
Call $L$ the complex line through $z':=\Phi(z)$ and $w':=\Phi(w)$.
Let $z_0 \in L \cap \partial \Omega$ be such that
$|z'-z_0| =\delta_{L\cap \Omega}(z').$ Let $P$ be the linear projection from $\C^n$ to $L$,
parallel to the complex tangent hyperplane to $\partial\Omega$ at $z_0$. Then $P(\Omega)$ is
a simply connected domain (see e.g.~\cite[Theorem 2.3.6]{APS}), and $z_0\in\partial P(\Omega)$.
Therefore,
\begin{multline*}
\tilde k_{D\cap U_0}(z,w) = \tilde k_{\Omega}(z',w') \ge \tilde k_{P(\Omega)}(z',w')
\\
\ge \frac14 \log \left( 1+ \frac{|z'-w'|}{\delta_{P(\Omega)}(z')} \right)
= \frac14 \log \left( 1+ \frac{|z'-w'|}{\delta_{L \cap \Omega}(z')} \right),
\end{multline*}
(for the second inequality see e.g. \cite[Proposition 3(ii)]{NT}).
By \cite[Propositions 4 and 6]{NPZ}, $\delta_{L \cap \Omega}(z') \lesssim \delta_\Omega(z')^{1/2m}$;
since $\Phi$ is biholomorphic in a neighborhood of $\overline{D\cap U_2}$, we have
$|z'-w'| \asymp|z-w|$ and $\delta_\Omega(z')=\delta_{D\cap U_0}(z)$,
so we finally obtain \eqref{mest} (the implicit constants are uniform over $D$ by a compactness argument).

{\bf Case 2.} $|z-w|^{2m} \ge \delta_D(z)$.

We may assume that $D\cap U_0$ is $\mathcal C^2$-smooth,
and that the projection $\pi$ to $\partial D$ is well defined on $U_0$.

We will follow the proof of \cite[Theorem 2.3]{FR}.
We need to bound from below the Kobayashi-Royden length of any path $\gamma$ such
that $\gamma(0)=z$ and $\gamma(1)=w$. If $\gamma([0,1]) \not\subset U_1$
(in particular if $w\notin U_1$), let
$t^*:= \min \{t \in [0,1]: \gamma(t) \notin U_1\}$. It will be enough to bound below
the length of $\gamma[0,t^*]$, so we can reduce ourselves to the case where $w\in \overline{U_1}$.

Let $\Phi$ be a holomorphic embedding such that $\Phi (D\cap U_0)=:\Omega$ is $\C$-convex.

Applying a result of K. Diederich and J.E. Fornaess about
supporting functions \cite{DF} to $\Omega$,
reducing $U_1$ as needed, we can find neighborhoods of $a$, $U_1 \Subset U_2 \Subset U_0$
such that for any $a' \in U_1$,
there exist $S_{\Phi(a')}$ holomorphic on $\C^n$, and $C, C'>0$ such that
\begin{equation}\label{peak}
-C'|\xi - \Phi(a')| \le \Re S_{ \Phi(a')} (\xi) \le\ - C |\xi - \Phi(a')|^{2m},
\end{equation}
\hfill $\xi \in \Phi(U_2), \mbox{ and } S_{ \Phi(a')}( \Phi(a'))=0.$
\smallskip\newline
We define a function $P_z$ holomorphic on $U_0$ by
\begin{equation}
\label{peakdef}
P_z (\zeta):= e^{S_{\Phi(\pi(z))}(\Phi(\zeta))}.
\end{equation}
Since $\Phi$ is a uniformly bilipschitz diffeormorphism on $U_2$ we then have, for $\zeta \in U_2$,
\begin{equation}
\label{hpeak}
\left| 1- P_z(\zeta)\right| \lesssim |\zeta-\pi(z)|
\mbox{ and }
1 - |P_z(\zeta)| \gtrsim |\zeta-\pi(z)|^{2m}.
\end{equation}
This means in particular that \cite[Lemma 2.2]{FR} can be applied, and it follows that
by \cite[Theorem 2.1]{FR} that there is $C_1>0$ such that for $z \in D\cap U_1$ and $X\in \C^n$,
$$
\kappa_{D\cap U_0} (z; X) \ge
\kappa_D (z; X) \ge (1-C_1 \delta_D (z)) \kappa_{D\cap U_0} (z; X).
$$
Therefore
\begin{equation}
\label{lengthgamma}
\int_0^1 \kappa_D (\gamma(t), \gamma'(t)) dt \ge
\int_0^1 \left( 1-C_1 \delta_D (\gamma(t))\right) \kappa_{D\cap U_0}  (\gamma(t), \gamma'(t)) dt.
\end{equation}
Let $\lambda := P_z \circ \gamma$. Then
$$
\kappa_{D\cap U_0}  (\gamma(t), \gamma'(t)) \ge
\kappa_{\D}  (\lambda(t), \lambda'(t)) \ge \frac{|\lambda'(t)|}{2(1-|\lambda(t))|}.
$$
On the other hand, by \eqref{hpeak},
\begin{multline*}
1-C_1 \delta_D (\gamma(t)) \ge 1 -C_1 |\gamma(t)-\pi(z)|
\\
\ge 1-C_1'\left( 1-|P_z(\gamma(t))| \right)^{1/2m}
=  1-C_1'\left( 1-|\lambda(t)| \right)^{1/2m}.
\end{multline*}
Collecting the estimates, the double right hand side in \eqref{lengthgamma} can be bounded below by
\begin{multline*}
\int_0^1 \frac{1-C_1'\left( 1-|\lambda(t)| \right)^{1/2m}}{1-|\lambda(t)|} |\lambda'(t)| dt
\ge\int_0^1 \frac{1}{1-|\lambda(t)|} \frac{d}{dt}|\lambda(t)| dt + O(1)
\\
=
\log \frac{1-|\lambda(1)|}{1-|\lambda(0)|} + O(1)
=
\log \frac{1-|P_z(w)|}{1-|P_z(z)|} + O(1).
\end{multline*}
By \eqref{hpeak}, $1-|P_z(z)|  \lesssim |z-\pi(z)| = \delta_D (z)$, while
$$
1-|P_z(w)|  \gtrsim |w-\pi(z)|^{2m} \ge \left(|w-z| - |z-\pi(z)| \right)^{2m}.
$$
Since $\delta_D (z) \le (C_0^{-1} |w-z|)^{2m} < \frac12 |w-z|$ for $C_0$ large enough,
we have $1-|P_z(w)| \gtrsim |w-z|^{2m}$ and the estimate we wanted is proved.
\end{proof*}

\section{Proof of Theorem \ref{greenlocal}}
\label{estgreen}

\noindent{\it Proof of Theorem \ref{greenlocal}, \eqref{pole}}.
Choose a bounded neighborhood $U_0$ of $a$ such that $D\cap U_0$ is
$\C$-convexifiable and $\mathcal C^2$-smooth.

{\bf Case 1.} $|z-w| \le\delta_D (w)^{1/2m}$.

We can choose a neighborhood $U\Subset U_0$ such that for any $w\in D\cap U$, then
$z\in D\cap U_0$ and $\delta_{D \cap U_0} (w) = \delta_{D} (w)$.

By Lemma \ref{tech}(ii), we have to prove that
$$
g_D (z,w) \ge \log \frac{|z-w|}{\delta_D(w)^{1/2m}} +O(1).
$$
We first reduce ourselves
to the study of $g_{D \cap U_0}$ by a standard argument.

\begin{lemma}
\label{greenloczn} Shrinking $U$ (if necessary), there is $C>0$
such that
\begin{equation}\label{grloc}
g_{D} (z,w)\ge g_{D \cap U_0} (z,w)-C,\quad z\in D \cap U_0,\ w\in D \cap U.
\end{equation}
\end{lemma}

Accepting this lemma, we apply Lempert's theorem to $D \cap U_0$
and obtain $g_{D} (z,w)\ge k_{D \cap U_0} (z,w)-C_a$.
By Theorem \ref{kobalocal},
$\tilde k_{D \cap U_0}  (z,w)$ satisfies \eqref{pole} (by shrinking $U$ once more if needed),
therefore
$$
k_{D \cap U_0} (z,w) \ge \log \frac{|z-w|}{\delta_D(w)^{1/2m}} +O(1),
$$
and we are done.

\noindent{\it Proof of Lemma \ref{greenloczn}.}

The proof is similar to that of \cite[Theorem 1]{Com}.

Let $\psi(z)=\log\frac{|z-a|}{\diam D}$ and
$U_1\Subset U_0\subsetneq D$ be a neighborhood of $a$ such that
and $\inf_{D\setminus U_0}\psi > c:=1+\sup_{D\cap\partial U_1}\psi.$ Fix $w\in
D\cap U_1$ and set
$$
d(w)=\inf_{z\in D\cap\partial U_1}g_{D\cap U_0}(z,w),\quad u(z,w)=
(c-\psi(z))d(w),\ z\in D.
$$
Since $u(z,w)\le g_{D\cap U_0}(z,w)$ for $z\in D\cap\partial U_1,$ and
$u(z,w)> 0>g_{D\cap U_0}(z,w)$ for $z\in \mathcal N \cap (D\cap U_0)$,
where $\mathcal N$ is a neighborhood of $\partial U_0$, the function
$$v(z,w)=\left\{\begin{array}{ll}
g_{D\cap U_0}(z,w),&w\in D\cap U_1\\
\max\{g_{D\cap U_0}(z,w),u(z,w)\},&w\in D\cap U_0\setminus U_1\\
u(z,w),&w\in D\setminus U_0
\end{array}\right.$$
is a plurisubharmonic function in $z$ with logarithmic pole at $w.$
Also $v(z,w)<cd(w),$ so $g_D(z,w)\ge v(z,w)-cd(w).$ Now \eqref{grloc} follows
by taking $U\Subset U_1$ and $C:=c\inf_{w\in D\cap U}d(w).$\quad\qed

{\bf Case 2.} $|z-w| \ge \delta_D (w)^{1/2m}$.

By Lemma \ref{tech}(i), we have to prove that
\begin{equation}
\label{boundfrac}
g_D (z,w)\gtrsim-\frac{\delta_D(w)}{|z-w|^{2m}}.
\end{equation}
By Theorem \ref{kobalocal} and Lempert's theorem,
\begin{equation}
\label{chen5}
g_{D \cap U_0} (z,w)\gtrsim-\frac{\delta_D(w)}{|z-w|^{2m}},\quad z\in D\cap U_0,\ w \in D\cap U.
\end{equation}
We will follow part of the proof of \cite[Lemma 3]{Ch}. The above inequality is analogous to
\cite[p. 29, inequality (5)]{Ch}.

Denote by $B(w,r)$ the ball with center $w$ and radius $r.$
Set $r_0:= \frac14 \mbox{dist}\, (U, D\setminus U_0)$, $\lambda:= \min\{r_0, |z-w|\}$,
so that
$$
D\cap B(w,\lambda) \subset D\cap B(w,2r_0) \subset D\cap U_0.
$$
Note that
\begin{equation}
\label{lambdaest}
\lambda \le |z-w| \le \frac{\diam D}{r_0} \lambda.
\end{equation}
Finally, let
$$
b:= - \inf \left\{ g_{D \cap U_0} (\zeta,w): |\zeta-w|=\lambda, \zeta \in D \right\} .
$$
Because of \eqref{chen5} and
\eqref{lambdaest},
\begin{equation}
\label{estb}
 b \lesssim \frac{\delta_D(w)}{\lambda^{2m}} \lesssim  \frac{\delta_D(w)}{|z-w|^{2m}}.
\end{equation}

Let
$$
v(\zeta) := b \frac{\log\frac{|\zeta-w|}{2r_0}}{\log\frac{2r_0}\lambda}.
$$
By construction, $v(\zeta) =0 > g_{D \cap U_0} (\zeta,w)$ when $\zeta \in D\cap \partial B(w,2r_0)$,
and $v(\zeta) =-b \le  g_{D \cap U_0} (\zeta,w)$ when $\zeta \in D\cap \partial B(w,\lambda)$.

Then we construct a plurisubharmonic function $u$ with logarithmic singularity at $w$ by setting
$$
u(\zeta):= \left\{
\begin{array}{l}
g_{D \cap U_0} (\zeta,w), \quad \zeta \in B(w,\lambda),
\\
\max\left\{ v(\zeta),  g_{D \cap U_0} (\zeta,w)\right\}, \quad \zeta \in B(w,2r_0) \setminus B(w,\lambda) ,
\\
v(\zeta), \quad \zeta \in D\setminus B(w,2r_0)
\end{array}
\right.
$$
By definition of $g_D$, $g_D \ge u - \sup_D u$.  We have
$$
\sup_D u \le \sup_D v \le b \frac{\log\frac{\diam D}{2r_0}}{\log\frac{2r_0}\lambda}
\le b \frac{\log\frac{\diam D}{2r_0}}{\log 2}  \lesssim  \frac{\delta_D(w)}{|z-w|^{2m}},
$$
by \eqref{estb}. On the other hand, if $\lambda = |z-w|$, then
$$
u(z) = g_{D \cap U_0} (z,w) \gtrsim -\frac{\delta_D(w)}{|z-w|^{2m}}
$$
by \eqref{chen5}, while if $\lambda = r_0 < |z-w|$, then
$$
u(z) \ge v(z) = b \frac{\log\frac{|z-w|}{2r_0}}{\log 2} \ge -b.
$$
Collecting the estimates, $g_D (z,w)  \gtrsim -\frac{\delta_D(w)}{|z-w|^{2m}}$.\quad\qed
\smallskip

\noindent{\it Proof of Theorem \ref{greenlocal}, \eqref{argu}.}
We choose $U_1$ small enough so that $\pi(z)$ is well defined whenever $z\in U_1$.

{\bf Case 1.} Suppose that $z\in U$ and $|z-w| \ge\delta_D(z)^{1/2m}$.

Shrinking $U_1,$ we may assume that $|z-w| \ge 8 \delta_D(z).$

We use the Diederich-Fornaess supporting functions \cite{DF} once again.
We take $U_1 \Subset U_2 \Subset U_0$ as before.
Reducing $U_1$ if needed, for any
$a' \in U_1 \cap \partial D$, there exist $S_{\Phi(a')}$ holomorphic on $\Omega$, and
$C,C'>0$ such that \eqref{peak} holds.

We set $\tilde \varphi_z (\zeta) := \Re S_{\Phi(\pi(z))}(\Phi(\zeta)) \in PSH_-(D\cap U_0)$.
Since $\Phi$ is a uniformly
bilipschitz diffeormorphism on $U_2$ we then have, for $\zeta \in U_2$,
\begin{equation}
\label{peakphi}
-C'|\zeta - a'| \le  \tilde \varphi_z(\zeta) \le - C |\zeta - a'|^{2m} \mbox{ and } \tilde \varphi_z(\pi(z))=0.
\end{equation}

We need to extend $\tilde \varphi_z$ to a global plurisubharmonic function on $D$.
We proceed as in \cite[p. 31]{Ch}. Let
$\eta := \sup_{z\in U_1} \sup_{\zeta \in \partial U_2} \tilde \varphi_z(\zeta) <0$.
We set $\varphi_z := \max ( \tilde \varphi_z, \eta/2)$ and extend it by $\eta/2$ on the whole of $D$.
Then $\varphi_z \in PSH_-(D)$ and satisfies the analogue of \eqref{peakphi}.

By the same argument as at the beginning of Case 2 of the proof of \eqref{pole},
the inequality we have to prove is the following analogue of \eqref{boundfrac}:
\begin{equation*}
g_D (z,w) \gtrsim - \frac{\delta_D(z)}{|z-w|^{2m}}.
\end{equation*}

\begin{lemma}
\label{bychen}
Let $w':= w + \frac{w-z}{|w-z|}$, $B_1:= B(w', 1 + |w-z|/2)$, $B_2:= B(w', 1 + 3|w-z|/4)$.
There is $c_0>0$
so that for any $w$, there exists $\rho_w \in \mathcal C^\infty (\C^n \setminus \{w\}, \R_-)$ with
logarithmic singularity at $w$, supported on $\overline B_2$, such that
$$
\partial \bar \partial \rho_w (\zeta) \ge - \frac{c_0}{|w-z|^2} \chi_{\overline B_2 \setminus B_1} (\zeta) \partial \bar \partial (|\zeta|^2).
$$
In particular, $\rho_w \in PSH (B_1 \cup (\C^n \setminus \overline B_2))$.
\end{lemma}
This lemma is proved in \cite[p. 31]{Ch}.

We construct a  function $\Phi$ with logarithmic pole at $w$ by setting
$$
\Phi (\zeta) := \frac{c_1}{|z-w|^{2m}} \left( \varphi_z (\zeta) + c_2 |\zeta - \pi(z) |^{2m}\right)
+ \rho_w (\zeta).
$$
By \eqref{peakphi} and because $D$ is bounded, we can choose $c_2>0$ such that $\Phi<0$ on $D$.

We want to choose $c_1>0$ so that $\Phi \in PSH(D)$.
We only need to check the case where $\zeta \in \overline B_2 \setminus B_1$. Then
$$
|\zeta - \pi(z)| \ge |\zeta - z| - \delta_D (z) \ge \frac14 |z-w| - \delta_D (z)  \ge \frac18 |z-w|.
$$
By the estimate on $\partial \bar \partial \rho_w $ from Lemma \ref{bychen}, the
fact that $\varphi_z \in PSH(D)$, and standard
computations,
\begin{multline*}
\partial \bar \partial \Phi (\zeta) \ge
\left(
 \frac{c_1}{|z-w|^{2m}} c_2 c_3 |\zeta - \pi(z) |^{2m-2} - \frac{c_0}{|w-z|^2}
\right) \partial \bar \partial |\zeta|^2
\\
\ge \left( \frac{c_1 c_2 c_3}{8^{2m-2}} -c_0
\right) \frac1{|w-z|^2} \partial \bar \partial |\zeta|^2,
\end{multline*}
where $c_3>0$ is a constant.  So we can choose $c_1>0$ to make this
form positive.  With these choices, $\Phi(\zeta) \le g_D(\zeta,w)$.

Since $\rho_w (z)=0$, using \eqref{peak} again,
$$
\Phi(z) = \frac{c_1}{|z-w|^{2m}} \left( \varphi_z (z) + c_2 \delta_D(z)^{2m}\right)
\ge - c_1 C'\frac{\delta_D(z) }{|z-w|^{2m}}.
$$

{\bf Case 2.} Suppose that $z\in B(a,r_1)$ and $|z-w| \le \delta_D(z)^{1/2m}$.

Then $|w-a| \le r_1 + r_1^{1/2m} =: r_2$.  Reducing $r_1$ if needed, we have
$B(a,r_2) \Subset U_0$, where
$U_0$ is a bounded neighborhood of $a$ such that
$D\cap U_0$ is $\C$-convexifiable and $\mathcal C^2$-smooth.
This implies that, by Lempert's theorem and \eqref{pole},
$$
\tilde g_{U_0 \cap D} (z,w) = \tilde g_{U_0 \cap D} (w,z) \ge
m \log \left( 1+C \frac{|z-w|}{\delta_D(z)^{1/2m}} \right).
$$
Since $\frac{|z-w|}{\delta_D(z)^{1/2m}} \le 1$,
by Lemma \ref{tech}(ii), this is equivalent to
$g_{U_0 \cap D} (z,w) \ge \log \frac{|z-w|}{\delta_D(z)^{1/2m}} +O(1)$.
By Lemma \ref{greenloczn}, the same estimate holds for $g_D (z,w)$,
and we are done for this case.\quad\qed
\smallskip

\section{Proofs of Theorems \ref{greenratio} and \ref{greengrowth}}
\label{proofgreen}

\noindent{\it Proof of Theorem \ref{greenratio}.}
Let
$$
\Delta_D (z,w):=\frac{|z-w|^2}{\delta_D(z)^{1/2m}\delta_D(w)^{1/2m}}.
$$
Using \eqref{sw}, it is enough to show that
$g_D(z,w)\gtrsim - \Delta_D (z,w)^{-2m}.$

Theorem \ref{greengrowth} implies that
$$\tilde g_D (z,w) \ge \log \left( 1+C'\Delta_D(z,w)\right)^m.$$

If $\Delta_D(z,w) \ge 1$, then $g_D(z,w)\gtrsim-\Delta_D (z,w)$ by Lemma \ref{tech}(i).

If $\Delta_D(z,w) \le 1$, then Lemma \ref{tech}(ii) implies that
$$g_D (z,w) \ge  \log \Delta_D(z,w)+O(1)\gtrsim-\Delta_D(z,w)^{-2m}.\qed$$

\noindent{\it Proof of Theorem \ref{greengrowth}.}
We follow an argument in \cite{Ca}, as adapted in \cite[Proof of Proposition 2]{Ch}.

The hypotheses of Theorem \ref{greenlocal} are met for any $a \in \partial D$.
By a compactness argument, this implies that there is $K \Subset D$ such that for
$z \in D\setminus K$, $w\in D$,
\begin{equation}
\label{ineqz}
\tilde g_D (z,w) \ge m \log \left(1 + C\frac{|z-w|}{\delta_D(z)^{1/2m}}\right).
\end{equation}
But when $z\in K$,  the right hand side of \eqref{ineqz} is bounded above
by $C' m C |z-w|$, while $\tilde g_D (z,w) \ge C'' |z-w|$, so $C$ can be chosen so that
\eqref{ineqz} holds for any $z , w\in D$. In the same way, changing $C$ again if needed,
we have for any $z , w\in D$,
\begin{equation}
\label{ineqw}
\tilde g_D (z,w) \ge m \log \left(1 + C\frac{|z-w|}{\delta_D(w)^{1/2m}}\right).
\end{equation}

If $|z-w|^{2m}\lesssim\max\{\delta_D(z),\delta_D(w)\}$,
then \eqref{growth} follows from \eqref{ineqz} and \eqref{ineqw}
by modifying the constant $C$.  Otherwise, by Lemma \ref{tech}(i), (2) is equivalent to
$$
g_D (z,w) \gtrsim - \frac{\delta_D(z)\delta_D(w)}{|z-w|^{4m}}.
$$
We may assume that $4\max\{\delta_D(z),\delta_D(w)\}\le|z-w|$. If $2|\zeta-\pi(z)|= |z-w|,$ then
$$
|\zeta-w| \ge |z-w| - |\zeta-\pi(z)| - |z-\pi(z)| \ge  \frac{|z-w|}{4}.
$$
Therefore, by \eqref{ineqw}, for those values of $\zeta$,
$g_D (\zeta,w) \gtrsim - \frac{\delta_D(w)}{|z-w|^{2m}}$.
For those same $\zeta$, the plurisubharmonic peak function $\varphi_z$ from the
proof of Theorem \ref{greenlocal}, \eqref{argu}, Case 1, verifies
$$
\varphi_z (\zeta) \le -C  |\zeta - \pi(z)|^{2m} = -C 2^{-2m} |z-w|^{2m},
$$
so,
$$
g_D (\zeta,w) \gtrsim - \frac{\delta_D(w)}{|z-w|^{4m}} \varphi_z (\zeta),
\quad\zeta \in D\cap \partial B(\pi(z),|z-w|/2).
$$
This inequality is trivially true on $\partial D,$ where $g_D(\zeta,w)=0$, and since
$g_D(\cdot,w)$ is a maximal plurisubharmonic function
on $D\setminus\{w\}$, it has to hold on $D\cap  B(\pi(z),|z-w|/2)$,
in particular at the point $z$, so
$$
g_D (z,w) \gtrsim - \frac{\delta_D(w)}{|z-w|^{4m}} \varphi_z (z)
\gtrsim -\frac{\delta_D(w)\delta_D(z)}{|z-w|^{4m}}.\quad\qed
$$

\end{document}